\documentclass[12pt,leqno]{article}
\usepackage{amsfonts}
\usepackage{amsmath, amsthm, amsfonts, amssymb, color}
\usepackage{mathrsfs}
\usepackage{color}
\theoremstyle{definition}

\newcommand{\scr}[1]{\mathscr #1}
\definecolor{wco}{rgb}{0.5,0.2,0.3}

\numberwithin{equation}{section} \theoremstyle{remark}

\newcommand{\ua}{\uparrow}

\title{
{\bf  Nash inequality for   Diffusion Processes Associated with Dirichlet Distributions}
\footnote{Supported in part by NNSFC (11771326,11726627).}
}
\author{{\bf     Feng-Yu Wang\footnote{Corresponding author: wangfy@tju.edu.cn}}\\
\footnotesize{Center for Applied Mathematics, Tianjin University, Tianjin 300072, China}\\
\footnotesize{Department of  Mathematics, Swansea University, Singleton Park, SA2 8PP, UK}\\
{\bf Weiwei Zhang}\\
\footnotesize{School of Mathematical Sciences, Beijing Normal University, Beijing 100875, China}
 }

\begin{document}
\def\tttext#1{{\normalfont\ttfamily#1}}
\def\R{\mathbb R}  \def\ff{\frac} \def\ss{\sqrt} \def\B{\mathbf B}
\def\N{\mathbb N} \def\kk{\kappa} \def\m{{\bf m}}
\def\dd{\delta} \def\DD{\Delta} \def\vv{\varepsilon} \def\rr{\rho}
\def\<{\langle} \def\>{\rangle} \def\GG{\Gamma} \def\gg{\gamma}
  \def\nn{\nabla} \def\pp{\partial} \def\EE{\scr E}
\def\d{\text{\rm{d}}} \def\bb{\beta} \def\aa{\alpha} \def\D{\scr D}
  \def\si{\sigma} \def\ess{\text{\rm{ess}}}
\def\beg{\begin} \def\beq{\begin{equation}}  \def\F{\scr F}
\def\Ric{\text{\rm{Ric}}} \def\Hess{\text{\rm{Hess}}}
\def\e{\text{\rm{e}}} \def\ua{\underline a} \def\OO{\Omega}  \def\oo{\omega}
 \def\tt{\tilde} \def\Ric{\text{\rm{Ric}}}
\def\cut{\text{\rm{cut}}} \def\P{\mathbb P}
\def\C{\scr C}     \def\E{\mathbb E}\def\y{{\bf y}}
\def\Z{\mathbb Z} \def\II{\mathbb I}
  \def\Q{\mathbb Q}  \def\LL{\Lambda}\def\L{\scr L}
  \def\B{\scr B}    \def\ll{\lambda} \def\a{{\bf a}} \def\b{{\bf b}}
\def\vp{\varphi}\def\H{\mathbb H}\def\ee{\mathbf e}\def\x{{\bf x}}
\def\gap{{\rm gap}}\def\PP{\scr P}\def\p{{\mathbf p}}\def\NN{\mathbb N}
\def\cA{\scr A} \def\cQ{\scr Q}\def\cK{\scr K}
\def\LS{C_{LS}}
\maketitle
\begin{abstract}  For any $N\ge 2$ and $\aa=(\aa_1,\cdots, \aa_{N+1})\in (0,\infty)^{N+1}$, let $\mu^{(N)}_{\aa}$ be the Dirichlet distribution with parameter $\aa$  on  the set
$\DD^{ (N)}:=  \{ x \in [0,1]^N:\  \sum_{1\le i\le N}x_i  \le 1 \}.$  The multivariate Dirichlet diffusion is associated with the     Dirichlet form
 $$\EE_\aa^{(N)}(f,f):=  \sum_{n=1}^N \int_{ \DD^{(N)}} \bigg(1-\sum_{1\le i\le N}x_i\bigg)  x_n(\pp_n f)^2(x)\,\mu^{(N)}_\aa(\d x)$$ with Domain $\D(\EE_\aa^{(N)})$ being the closure of $C^1(\DD^{(N)})$.  We prove the Nash   inequality
 $$\mu_\aa^{(N)}(f^2)\le C  \EE_\aa^{(N)}(f,f)^{\ff p{p+1} }\mu_\aa^{(N)} (|f|)^{\ff 2 {p+1}},\ \ f\in \D(\EE_\aa^{(N)}), \mu_\aa^{(N)}(f)=0$$ for some constant $C>0$ and
 $p= (\aa_{N+1}-1)^+ +\sum_{i=1}^N 1\lor (2\aa_i),$ where the constant $p$  is sharp when  $\max_{1\le i\le N} \aa_i  \le \ff 1 2$ and $\aa_{N+1}\ge 1$.   This Nash inequality also holds   for the corresponding  Fleming-Viot process.
 \end{abstract} \noindent
   \vskip 2cm

\section{Introduction}

 Let $ N\ge 1$ be a natural number, and let  $ \aa=(\aa_1,\cdots, \aa_{N+1})\in (0,\infty)^{N+1}.$   The Dirichlet distribution $\mu^{(N)}_{\aa}$ with parameter $\aa$ is a probability measure on the set
 $$\DD^{(N)}:= \bigg\{ x=(x_i)_{1\le i\le N}\in [0,1]^N:\ |x|_1:=\sum_{i=1}^N x_i\le 1\bigg\}$$ with density function
\beq\label{RR} \rr(x):= \ff{\GG(|\aa|_1)}{\prod_{1\le i\le N+1} \GG(\aa_i)} (1-|x|_1)^{\aa_{N+1}-1}\prod_{1\le i\le N} x_i^{\aa_i-1},\ \ x=(x_i)_{1\le i\le N}\in \DD^{(N)},\end{equation}
where $|\aa|_1:=\sum_{i=1}^{N+1} \aa_i$. 
This   distribution  arises naturally in Bayesian inference as conjugate prior  for categorical distribution, and it describes   the distribution of  allelic frequencies in  population genetics,  see for instance \cite{ConMoi69,P1,P2}.

 To investigate stochastic dynamics converging to $\mu^{(N)}_{\aa}$, different models of diffusion processes have been proposed. In this paper, we investigate functional inequalities of these diffusions.

 In the following three subsections,   we   first   briefly recall some facts on functional inequalities for Dirichlet forms, as well as   known results for    diffusion processes associated with the Dirichlet distribution, then   propose    problems  in the direction and  state the main result of the paper.

  \subsection{Functional inequalities}
 In general, let $(\EE,\D(\EE))$ be a conservative symmetric Dirichlet form on $L^2(\mu)$ for some probability space $(E,\F,\mu)$,   let $(L, \D(L))$ be the associated Dirichlet operator, and let $P_t:=\e^{tL}, t\ge 0, $ be the Markov semigroup.  The following is a brief summary from \cite{Wbook} for the Poincar\'e, log-Sobolev, super Poincar\'e and Nash inequalities, see also \cite{Bakry, Gross} and references within.

 Firstly, we consider the spectral gap     of $L$: $\gap(L)$ is the largest constant $C>0$ such that the Poincar\'e inequality
 \beq\label{P0}\mu(f^2)\le \ff 1 C \EE(f,f),\ \ f\in \D(\EE), \mu(f)=0\end{equation} holds.  In case this inequality is not available,  we say that $L$ does not have spectral gap, and denote $\gap(L)=0$. The Poincar\'e inequality \eqref{P0} is equivalent to the $L^2$-exponential  convergence of $P_t$: $$\|P_tf\|_{L^2(\mu)}\le \e^{-Ct} \|f\|_{L^2(\mu)},\ \  t\ge 0, f\in L^2(\mu), \mu(f)=0.$$

 Next,  we consider the log-Sobolev constant $\LS(L)$, which is the largest positive constant $C$ such that  the log-Sobolev inequality
\beq\label{LS0} \mu(f^2\log f^2)\le \ff 2 C \EE(f,f),\ \ f\in \D(\EE), \mu(f^2)=1 \end{equation} holds. We have $\LS(L)\le \gap(L)$. In general,  \eqref{LS0} implies the exponential decay of $P_t$ in entropy:
$$ \mu((P_t f)\log P_t f)\le \e^{-Ct }\mu(f\log f),\ \ t\ge 0, f\in \B^+(E), \mu(f)=1, $$ and in the diffusion setting they are equivalent. Moreover, the log-Sobolev inequality
\eqref{LS0} holds for some constant $C>0$ if and only if $P_t$ is hypercontractive, i.e. $\|P_t\|_{L^2(\mu)\to L^4(\mu)}=1$ for large enough $t$.

Finally, we say that $(\EE,\mu)$ satisfies the super Poincar\'e inequality with rate function $\bb: (0,\infty)\to (0,\infty),$ if
\beq\label{SP0} \mu(f^2)\le r\EE(f,f) +\bb(r) \mu(|f|)^2,\ \ r>0, f\in \D(\EE).\end{equation} This inequality is equivalent to the uniform integrability of $P_t$, i.e. $P_t$ has zero tail norm:
$$\|P_t\|_{tail}:=\lim_{R\to\infty} \sup_{\mu(f^2)\le 1} \mu((P_tf)^2 1_{\{|P_tf|\ge R\}}) =0,\
 \  t>0.$$ When $P_t$ has a heat kernel with respect to $\mu$, it is also equivalent to the absence of the essential spectrum of $L$ (i.e. the spectrum of  $L$ is purely discrete).
 The super Poincar\'e inequality generalizes the classical Sobolev/Nash type inequalities. For instance, when $\gap(L)>0$, \eqref{SP0} with $\bb(r)= \e^{c(1+r^{-1})}$ for some $c>0$ is
 equivalent to the log-Sobolev inequality \eqref{LS0} for some constant $C>0$; while for a constant $p>0$, \eqref{SP0} with $\bb(r)=  c(1+r^{-p})$ holds for some $c>0$ if and only if
 the Nash inequality
 \beq\label{NS}  \mu (f^2)\le C  \EE (f,f)^{\ff p{p+1} }\mu  (|f|)^{\ff 2 {p+1}},\ \ f\in \D(\EE ), \mu (f)=0\end{equation}
 holds for some constant $C>0$, they are also equivalent to
 $$\|P_t-\mu\|_{L^1(\mu)\to L^\infty(\mu)} \le \ff{c'}{(t\land 1)^p} \e^{-\gap(L)t},\ \ t>0.$$ The later implies the hypercontractivity of $P_t$, and hence the log-Sobolev inequality \eqref{LS0} for some constant $C>0$. 

 \subsection{Diffusion processes  associated with Dirichlet distributions} In this part, we  recall existing results on functional inequalities for some  diffusion processes on $\DD^{(N)}$, which are reversible with respect to the Dirichlet distribution  $\mu_\aa^{(N)}$.

 When $N=1$, the most popular model is the Wright-Fisher diffusion on the interval $[0,1]$ generated by
 $$L_\aa^{(1)} := x(1-x) \pp_x^2 +\{\aa_1(1-x) -\aa_2x\} \pp_x.$$ The associated Dirichlet form is the closure of $(\EE_\aa^{(1)}, C^1([0,1]))$ given by
 $$\EE_\aa^{(1)}(f,g)= \int_0^1 x(1-x) f'(x)g'(x) \mu_\aa^{(1)}(\d x),\ \ f,g\in C^1([0,1]).$$
Due to \cite{SN}, we have $\gap(L_\aa^{(1)}) =\aa_1+\aa_2$, and \cite[Lemma 2.7]{S} shows that
 $\LS(L)\ge \ff{\aa_1\land \aa_2}{160}.$ Moreover, according to \cite[Theorem 2.2]{FW14}, $(\EE_\aa^{(1)},\mu_\aa^{(1)})$ satisfies the super Poincar\'e inequality with $\bb(r)=c(1+r^{-(\ff 1 2 \lor \aa_1\lor \aa_2)})$ for some constant $c>0$, and hence, the Nash inequality holds for $p= \ff 1 2\lor \aa_1\lor \aa_2$, which is sharp in the sense that the super Poincar\'e inequality does not hold if
 $\lim_{r\to 0} \bb(r) r^{\ff 1 2 \lor \aa_1\lor \aa_2}=0.$

 When $N\ge 2$, we consider the following three different generalizations of the Wright-Fisher diffusion arising from population genetics, see e.g. \cite{FW07,FW14,Mi1,Mi2,S}.

\paragraph{A. Fleming-Viot process.}    Let $|\aa|_1=\sum_{i=1}^{N+1}\aa_i$ and denote $\pp_i=\pp_{x_i}, 1\le i\le N$.  Consider the diffusion process on $\DD^{(N)}$  generated by
$$L_{FV}^{\aa,N} := \sum_{i,j=1}^N x_i(\dd_{ij}-x_j) \pp_i\pp_j+ \sum_{i=1}^N (\aa_i-|\aa|_1x_i)\pp_i, $$
where $\dd_{ij}= 1$ if $i=j$; $=0$ otherwise.
The associated Dirichlet form is the closure of $(\EE_{FV}^{\aa,N},C^1(\DD^{(N)}))$ given by
$$\EE_{FV}^{\aa,N}(f,g)= \int_{\DD^{(N)}} \sum_{i,j=1}^N x_i(\dd_{ij}-x_j)\{(\pp_i f)(\pp_j g)\}(x) \mu_\aa^{(N)}(\d x),\ \ f,g\in C^1(\DD^{(N)}).$$
Again due to \cite{SN, S} we have
$$\gap(L_{FV}^{\aa,N})=|\aa|_1, \ \  \LS(L)\ge \ff{1}{160}\min_{1\le i\le N+1} \aa_i.$$
  However,  {\bf the Nash inequality is unknown.}

\paragraph{B. GEM process.}  Let $\bb_i= \sum_{j=i+1}^{N+1} \aa_j,\ 1\le i\le N$. Then $\mu_\aa^{(N)}= \Pi_{\aa,\bb}$, the GEM distribution with parameter $(\aa,\bb)$, see e.g. \cite{FW07}.  For $x\in \DD^{(N)}$ and $1\le i,j\le N$, let
\beg{align*}& a_{ij}(x)= x_ix_j \sum_{k=1}^{i\land j} \ff{\big\{\dd_{ki} \big(1-\sum_{1\le l\le k-1} x_l\big)-x_k\big\}\cdot \big\{\dd_{kj} \big(1-\sum_{1\le l\le k-1} x_l\big)-x_k\big\}}{x_k\big(1-\sum_{1\le l\le k}x_l\big)},\\
& b_i(x)= x_i \sum_{k=1}^i  \ff{\big\{\dd_{ki} \big(1-\sum_{1\le l\le k-1} x_l\big)-x_k\big\}\cdot \big\{\aa_k \big(1-\sum_{1\le l\le k-1} x_l\big)-\bb_ix_k\big\}}{x_k\big(1-\sum_{1\le l\le k}x_l\big)}.\end{align*}
The corresponding GEM process introduced in \cite{FW07} is the diffusion process on $\DD^{(N)}$ generated by
$$L_{GEM}^{\aa,N}:= \sum_{i,j=1}^N a_{ij}\pp_i\pp_j + \sum_{i=1}^N b_i \pp_i,$$ and the associated Dirichlet form is the closure of   $(\EE_{GEM}^{\aa,N},C^1(\DD^{(N)}))$:
$$\EE_{GEM}^{\aa,N}(f,g)= \int_{\DD^{(N)}} \sum_{i,j=1}^N a_{ij}(x) \{(\pp_i f)(\pp_j g)\}(x) \mu_\aa^{(N)}(\d x),\ \ f,g\in C^1(\DD^{(N)}).$$
According to \cite[Theorem 3.1]{FW07}, we have
$$\gap(L_{GEM}^{\aa,N})=\aa_N+\aa_{N+1},\ \ C_{LS}(L_{GEM}^{\aa,N}) \ge \ff{1}{160}\min_{1\le i\le N+1} \aa_i.$$
Moreover, applying \cite[(1.4)]{FW14} for $a_i=\aa_i, b_i=\bb_i:=\sum_{j=i+1}^{N+1}\aa_j$, and using \cite[(2.24)]{FW14}, we see that the heat kernel $p_t^{GEM}(x,y)$ of the present GEM process with respect to $\mu_\aa^{(N)}$ satisfies
$$c_1 t^{-\sum_{i=1}^N(\ff 1 2\lor \aa_i\lor \bb_i)} \le \sup_{x,y\in \DD^{(N)}} p_t^{GEM}(x,y)\le c_2 t^{-\sum_{i=1}^N(\ff 1 2\lor \aa_i\lor \bb_i)},\ \ t\in (0,1]$$ for some constants $c_2>c_1>0$. So,
there exists a constant $C>0$ such that the Nash inequality \eqref{NS}  holds for $(\EE_{GEM}^{\aa,N}, \mu_\aa^{(N)})$ replacing $(\EE, \mu)$   with
$$p=\sum_{i=1}^N \max\Big\{ \ff 1 2, \aa_i, \sum_{ i+1\le j\le N+1}\aa_j\Big\},$$ which is sharp in the sense that the Nash inequality fails when this $p$ is replaced by any smaller constant.

\paragraph{C. Multivariate Dirichlet diffusion.}  This process was introduced   in  \cite{Jac01}, and was used in \cite{BR} to describe  a fluctuating
ensemble of $N$ variables subject to a conservation principle.   It can be constructed as the unique solution to  the following SDE on $\DD^{(N)}$:
\beq\label{E1} \d X_i(t)= \big\{\aa_i (1-|X(t)|_1)- \aa_{N+1}X_i(t)\big\}\d t + \ss{2 (1-|X(t)|_1)X_i(t)}\,\d B_i(t),\  \ 1\le i\le N,\end{equation}
where   $B(t):=(B_1(t),\cdots, B_N(t))$ is the $N$-dimensional Brownian motion. The infinitesimal generator of the diffusion is
 $$L^{(N)}_\aa:= \sum_{1\le n\le N} \Big(x_n (1-|x|_1)  \pp_n^2 + \big\{\aa_n (1-|x|_{1})-\aa_{N+1}x_n\big\}\pp_n\Big), $$
and the associated Dirichelt form  is the  closure of $(\EE_\aa^{(N)}, C^1(\DD^{(N)}))$:
$$\EE_{\aa}^{(N)}(f,g)= \int_{\DD^{(N)}} \sum_{i=1}^N x_i(1-|x|_1) \{(\pp_i f)(\pp_j g)\}(x) \mu_\aa^{(N)}(\d x).$$
According to \cite[Theorem 1.1]{FMW}, we have
$$\gap(L_\aa^N)=\aa_{N+1}.$$ Not that when $N=1$,   $\gap(L_\aa^{(1)})= \aa_1+\aa_2>\aa_2$.

Moreover, the whole spectrum of
$L_\aa^{(N)}$ has been characterized in \cite{FMW}. In particular, the essential spectrum is empty, so that the super Poincar\'e inequality
\beq\label{SP} \mu_\aa^{(N)} (f^2)\le r\EE_\aa^{(N)} (f,f) +\bb(r) \mu_\aa^{(N)}(|f|)^2,\ \ r>0, f\in C^1(\DD^{(N)})\end{equation} holds  for
some   function $\bb: (0,\infty)\to (0,\infty).$
However,  there is no any estimates on $\bb(r)$ and hence, both {\bf the log-Soblev and  the Nash inequalities are unknown}.

\subsection{Questions and the Main result}

According to  the last subsection,      the following two things remain unknown.

\beg{enumerate} \item[$(Q_1)$] Nash inequality for the Fleming-Viot and  multivariate Dirichlet diffusion processes.
\item[$(Q_2)$] Estimates  on  the log-Sobolev constant for the multivariate Dirichlet diffusion, and   the sharp log-Sobolev constant for the Wright-Fisher/Fleming-Viot/GEM processes.   \end{enumerate}
In this paper, we only investigate $(Q_1)$, and the main result is the following.

\beg{thm}\label{T1.1} Let $N\ge 2$.  \beg{enumerate}
\item[$(1)$] There exists a constant $C>0$ such that the Nash inequality
   \beq\label{SB} \mu_\aa^{(N)}(f^2)\le C  \EE_\aa^{(N)}(f,f)^{\ff p{p+1} }\mu_\aa^{(N)} (|f|)^{\ff 2 {p+1}},\ \ f\in \D(\EE_\aa^{(N)}), \mu_\aa^{(N)}(f)=0\end{equation}  holds for
   $p=p_\aa:=   \sum_{i=1}^N 1\lor (2\aa_i)+(\aa_{N+1}-1)^+,$  and the inequality remains true for $\EE_{FV}^{(\aa,N)}$ replacing $\EE_\aa^{(N)}$.
\item[$(2)$] If $\eqref{SB}$ holds for some constant $C>0$, then
$$p\ge \tt p_\aa:= \max\Big\{ \max_{1\le i\le N+1} \sum_{j\ne i, 1\le j\le N+1} \aa_j,\   \aa_{N+1}  + \max_{1\le i\le N} \sum_{j\ne i, 1\le j\le N} (1\lor \aa_j)\Big\}.$$
\item[$(3)$] If $\eqref{SB}$ with $\EE_{FV}^{\aa,N}$ replacing $\EE_\aa^{(N)}$ holds for some constant $C>0$, then
$$p\ge p_\aa':=  \max\Big\{   \sum_{ 1\le j\le N} \aa_i,\   \ff 1 2 \aa_{N+1}  + \ff 1 2\max_{1\le i\le N} \sum_{j\ne i, 1\le j\le N} (1\lor \aa_j)\Big\}.$$
   \end{enumerate}   \end{thm}

   \paragraph{Remark 1.2.} (1) Let $p_c$ be the smallest positive constant $p>0$ such that \eqref{SB} holds for some constant $C>0$, then assertions (1)-(2) in Theorem \ref{T1.1} imply
   $p_c\in [\tt p_\aa,p_\aa]$. In particular, when $\max_{1\le i\le N} \aa_i  \le \ff 1 2$ and $\aa_{N+1}\ge 1$, we have
$p_c= N+\aa_{N+1}-1$; that is, in this case the Nash inequality presented in Theorem \ref{T1.1}(1) is sharp for $\EE_\aa^{(N)}$. But the sharpness for
$\EE_{FV}^{\aa,N}$ is unknown. 

(2) As mentioned in the end of Subsection 1.1 that the Nash inequality \eqref{SB} implies that the log-Sobolev inequality
$$\mu_\aa^{(N)}(f^2\log f^2)\le C \EE_\aa^{(N)}(f,f),\ \ \ f\in \D(\EE_\aa^{(N)}), \mu_\aa^{(N)}(f^2)=1$$ holds for some constant $C>0$. However, in the moment we do not have any explicit estimate on the log-Sobolev constant $C_{LS}(L_\aa^{(N)}).$ 

   (3) Consider  the infinite-dimensional setting where $N=\infty$ and $\aa=(\aa_i)_{i\in \bar {\mathbb N}}$ with $|\aa|_1:= \sum_{i\in \bar {\mathbb N}} \aa_i<\infty$.
   According to \cite{S, FMW}, we  have  $$\gap(L_{FV}^{\aa,\infty})=  |\aa|_1,\ \ \gap(L_\aa^{(\infty)}) =\aa_\infty.$$
Next,  \cite[Theorem 3.5]{S} shows that  the set
$$D_0:=\big\{f\in \D(\EE_{FV}^{\aa,\infty}):\ \mu_\aa^{(\infty)}(f^2)+  \EE_{FV}^{\aa,\infty}(f,f)\le 1\big\}$$ is not uniform integrable in $L^2(\mu_\aa^{(\infty)})$, so that the super Poincar\'e inequality is not available for $(\EE_{FV}^{\aa,\infty},\mu_\aa^{(\infty)})$. Indeed, by  \cite[Theorem 1.2]{W00b} (see also \cite{W00a, Wbook}), if there exists $\bb: (0,\infty)\to (0,\infty)$ such that
$$ \mu_\aa^{(\infty)} (f^2)\le r \EE_{FV}^{\aa,\infty}(f,f) +\bb(r) \mu_\aa^{(\infty)} (|f|)^2,\ \ r>0, f\in \D(\EE_{FV}^{\aa,\infty}),$$
  then there exists
a positive increasing function $F$ on $[0,\infty)$ with $F(r)\uparrow\infty$ as $r\uparrow\infty$ such that
$$\mu_\aa^{(\infty)} (f^2F(f^2))\le  \EE_{FV}^{\aa,\infty}(f,f),\ \ f\in \D( \EE_{FV}^{\aa,\infty}),\ \mu_\aa^{(\infty)}(f^2)\le 1,$$ and hence $D_0$ is uniform integrable in $L^2(\mu_\aa^{(\infty)})$.    Since $ \EE_{FV}^{\aa,\infty}\ge \EE_\aa^{(\infty)},$ see the beginning of Section 3 for finite $N$, the super Poincar\'e inequality  is   invalid for $\EE_\aa^{(\infty)}$  neither.

\

 To prove Theorem \ref{T1.1}, we will present a localization theorem in Section 2, which enables one to establish the super Poincar\'e inequality \eqref{SP} by using local inequalities. A complete proof of Theorem \ref{T1.1} will be addressed in Section 3 and Section 4.

 \section{Preparations }

 To establish \eqref{SP} with  an explicit rate function $\bb$, the main difficulty   comes from the singularity of the density $\rr(x)$ as well as the degeneracy of the diffusion coefficient on the boundary $$\pp \DD^{(N)}=\Big\{x=(x_i)_{1\le i\le N}\in \DD^{(N)}:\ \min\{x_i: 1\le i\le N+1\}=0\Big\},\ \ x_{N+1}:=1-\sum_{i=1}^N x_i.$$  To overcome such  difficulties, a localization result has been  presented in  \cite[Theorem 3.4.6]{Wbook}. However, this result is less sharp and inconvenient for application to the present model. So, in this section we   give   a new version of this result.    We will also present an additivity property of the super Poincar\'e inequality, which is more or less trivial but will be used to establish local super Poincar\'e inequalities in the proof of Theorem \ref{T1.1}(1).

\subsection{A localization result}

 Let $(E,\F,\mu)$ be a separable complete probability space, and let
$(\EE,\D(\EE))$ be a conservative symmetric local Dirichlet form on
$L^2(\mu)$ as the closure of 
$$\EE(f,g)=\mu(\GG(f,g)),\ \ f,g\in \D_0(\GG),$$ 
where $\GG: \D(\GG)\times \D(\GG)\to \scr B(E)$ is a
  positive definite symmetric bilinear mapping,   $\scr B(E)$ is the set of all
$\mu$-a.e. finite measurable real functions on $E$,   $\D(\GG)$ is a sub-algebra of   $\B(E)$, and $ \D_0(\GG):=\{f\in \D(\GG): f^2,\GG(f,f)\in L^1(\mu)\}$ such that
\beg{enumerate} \item[(a)]  $\D_0(\GG)$ is dense in $L^2(\mu)$.
\item[(b)] $\D(\GG)$ is closed under combinations with   $\psi\in C([-\infty,\infty])$ such that $\psi$ is $C^1$ in $\R$ and  $\psi'$ has compact support, and $\GG(\psi\circ f, g)= \psi'(f) \GG(f,g)\ \mu$-a.e.  for $f,g\in \D(\GG)$.
\item[(c)] $\GG(fg,h)= g\GG(f,h)+ f\GG(g,h)\ \mu$-a.e. for $f,g,h\in \D(\GG)$.  \end{enumerate} We aim to establish the super Poincar\'e inequality
\eqref{SP0} with an  explicit $\bb: (0,\infty)\to (0,\infty).$

  \beg{thm}\label{T2.1} Let $\phi\in \D(\GG)$ be an unbounded  nonnegative   function such that
   \beq\label{HU} h(s):= \sup_{D_s} \GG(\phi,\phi)<\infty,\ \ D_s:=\{\phi\le s\},\ \ s\ge 0,\end{equation}   where $\sup_{\emptyset}=0$ by convention.  If there exists $s_0\ge 1$ such that for every $s\ge s_0$, the local super Poincar\'e inequality
  \beq\label{LSP2} \mu(f^2)\le r\EE(f,f)+ \bb_s(r) \mu(|f|)^2,\ \ r>0, f\in \D(\EE), f|_{D_s^c}=0 \end{equation} holds for some  decreasing function  $\bb_s: (0,\infty)\to (0,\infty)$,
    and for $s\ge s_0$
    \beq\label{D2} 0<\ll(s):= \inf\{\EE(f,f): \mu(f^2)=1, f|_{D_s}=0\}\uparrow \infty\ \text{as}\ s\uparrow\infty.\end{equation}
    Then
   \beq\label{SR}s_r:= \inf\{s\ge s_0: \ll(s)\ge 8r^{-1}\}\in (0,\infty),\ \ r>0,\end{equation} and there exists a constant $c>0$ such that the super Poincar\'e inequality
    $\eqref{SP0}$ holds for
  \beq\label{BT} \bb(r):= c+ \Big(2 +\ff{rh(2s_r)}{s_r^2}\Big)\bb_{3s_r}\Big(\ff{r}{8+ 2 rh(2s_r)s_r^{-2}}\Big),\ \ r>0.\end{equation}
  \end{thm}
\beg{proof} By condition (a), it suffices to consider  $f\in \D_0(\GG)$. For any $s\ge s_0$ and small $\vv\in (0,1)$, let
$\varphi_i\in C^1([0,\infty])$ with $0\le \varphi_i\le 1, |\varphi_i'(s)|\le (1+\vv) s^{-1}, i=1,2$ such that
$$\varphi_1|_{[0, s]}=0,\ \ \varphi_1|_{[2s,\infty]}=1;\ \ \varphi_2|_{[0,2s]}=1,\ \ \varphi_2|_{[3s, \infty]}=0.$$
Let $f_i= f \cdot \varphi_i\circ\phi, 1\le i\le 2.$ Then
$f^2\le f_1^2+f_2^2$ and by conditions (b) and (c),
\beg{align*} \GG(f_1,f_1)&\le 2\GG(f,f)+ 2(1+\vv)^2 f^2 s^{-2} 1_{\{\phi\le 2 s\}} h(2s)\\
& \le  2\GG(f,f) + \ff{2(1+\vv)^2h(2s)}{s^2} f_2^2,\\
 \GG(f_2,f_2)&\le 2 \GG(f,f) + \ff 2 {s^2} (1+\vv)^2  f^2.\end{align*} In particular, $f_1,f_2\in \D_0(\GG)\subset \D(\EE)$.
Combining these with \eqref{LSP2} and \eqref{D2}, we obtain
\beq\label{D*} \beg{split}& \mu(f^2)  \le \mu(f_1^2)+ \mu(f_2^2) \le \ff 2 {\ll(s)} \EE(f,f) +\Big(1+\ff{2(1+\vv)^2h(2s)}{s^2\ll(s)}\Big)\mu(f_2^2)\\
&\le \Big\{\ff 2 {\ll(s)}+ 2t \Big(1+\ff{2(1+\vv)^2h(2s)}{s^2\ll(s)}\Big)\Big\} \EE(f,f)+ \ff{2(1+\vv)^2 t}{s^2} \Big(1+\ff{2(1+\vv)^2h(2s)}{s^2\ll(s)}\Big)\mu(f^2)\\
&\qquad + \Big(1+\ff{2(1+\vv)^2h(2s)}{s^2\ll(s)}\Big)\bb_{3s}(t) \mu(|f|)^2,\ \ t>0.\end{split}\end{equation} Let $r\in (0,1]$, and let $s_r$ be in \eqref{SR}. We have $\ll(s_r)\ge 8r^{-1},$ so that
$$\ff2{\ll(s_r)}\le \ff r 4,\ \ t_r:= \ff{r}{8+ 16 h(2s_r)/[s_r^2 \ll(s_r)]}\ge \ff{r}{8+ 2rs_r^{-2} h(2s_r)},$$ and
  when $\vv>0$ is small enough,
\beg{align*}& (1+\vv)^2\le 2,\ \
  \ff{2(1+\vv)^2h(2s_r)}{s_r^2\ll(s_r)}\le  \ff{rh(2s_r)}{2 s_r^2},\\
&2(1+\vv)^2t_r\Big(1+\ff{2(1+\vv)^2h(2s_r)}{s_r^2\ll(s_r)}\Big) \le \ff {3r} 8,\\
& \ff{2(1+\vv)^2 t_r}{s_r^2} \Big(1+\ff{2(1+\vv)^2h(2s_r)}{s_r^2\ll(s_r)}\Big) \le \ff{3 r} 8 \le \ff 3 8,\ \ r\in (0,1].\end{align*}
Combining these with \eqref{D*} we arrive at
$$\mu(f^2)\le \ff {5r} 8 \EE(f,f) +\ff 3 8 \mu(f^2)+ \Big(1+ \ff{rh(2s_r)}{2s_r^2}\Big) \bb_{3 s} \Big(\ff{r  }{8  + 2 r s_r^{-2}h(2 s_r)}\Big) \mu(|f|)^2,\ \ r\in (0,1].$$ Therefore,
$$\mu(f^2)\le r\EE(f,f) +  \Big(2+ \ff{rh(2s_r)}{s_r^2}\Big) \bb_{3 s} \Big(\ff{r  }{8   + 2 r s_r^{-2}h(2 s_r)}\Big) \mu(|f|)^2,\ \ r\in (0,1].$$
Since for the super Poincar\'e inequality we may take decreasing $\bb$, this finishes the proof.

\end{proof}

\subsection{Additivity of super Poincar\'e inequality}
For every $1\le i\le N$, let $(\EE_i,\D(\EE_i))$ be a symmetric Dirichlet form  on $L^2(\mu_i)$ over a $\si$-finite measure space $(E_i,\F_i,\mu_i)$. Let $\mu=\prod_{i=1}^N \mu_i$, and let $\D(\EE)$ be the class of $f\in L^2(\mu)$ such that for any $1\le i\le N$ and $\big(\prod_{j\ne i} \mu_j\big)$-a.e. $x$, we have $f(x,\cdot)\in \D(\EE_i)$ and 
$$\EE(f,f):= \sum_{i=1}^N \int_{\prod_{j\ne i} E_j} \EE_i(f(x,\cdot),f(x,\cdot))  \Big(\prod_{j\ne i}\mu_j\Big)(\d x)<\infty.$$
  Consider the following Dirichlet form on $L^2(\mu)$:
$$\EE(f,g):= \sum_{i=1}^N \int_{\prod_{j\ne i} E_j} \EE_i(f(x,\cdot),g(x,\cdot))  \Big(\prod_{j\ne i}\mu_j\Big)(\d x),\ \ f,g\in \D(\EE).$$
The following additivity property is a simple consequence of the equivalence between the heat kernel upper bound and the super Poincar\'e  inequality. 

\beg{prp}\label{PAD} Let $\{p_i\}_{1\le i\le N}\subset (0,\infty)$ such that for any $1\le i\le N$, the super Poincar\'e inequality
\beq\label{SPI} \mu_i(f^2)\le r \EE_i(f,f) + c_i (1+r^{-p_i}) \mu_i(|f|)^2,\ \ r>0, f\in \D(\EE_i)\end{equation} holds for some constant $c_i>0$. Then
there exists a constant $c>0$ such that
\beq\label{SPN} \mu(f^2)\le r \EE(f,f) + c \big(1+r^{-\sum_{i=1}^Np_i}\big) \mu(|f|)^2,\ \ r>0, f\in \D(\EE).\end{equation}
\end{prp}

\beg{proof} Let $P_t^i$ be the (sub) Markov semigroup associated with $(\EE_i,\D(\EE_i))$. By \cite[Theorem 3.3.15(2)]{Wbook}, \eqref{SPI} implies that $P_t^i$ has a density $p_t^i(x_i,y_i)$ with respect to $\mu_i$ such that
$${\rm ess}_{\mu_i\times \mu_i}\sup p_t^i\le C_i (1+t^{-p_i}),\ \ t>0$$ holds for some constant $C_i>0$. Then the semigroup $P_t$ associated with $(\EE,\D(\EE))$ has the density
$$p_t(x,y):=\prod_{i=1}^N p_t^i(x_i,y_i),\ \  x=(x_1,\cdots, x_N), y=(x_1,\cdots, y_N)\in E:=\prod_{i=1}^N E_i$$  with respect to $\mu$, and
$${\rm ess}_{\mu\times \mu}\sup p_t\le C \big(1+t^{-\sum_{i=1}^Np_i}\big),\ \ t>0$$ holds for some constant $C>0$.
By \cite[Theorem 3.3.15(2)]{Wbook} again, this implies \eqref{SPN} for some constant $c>0$. \end{proof}

\section{Proof of Theorem \ref{T1.1}(1)}

We first observe that  for any $f\in C^1(\DD^{(N)})$,
\beg{align*} &\sum_{i,j=1}^N x_i(\dd_{ij} -x_j) (\pp_if)\pp_j f = \sum_{i=1}^N x_i (\pp_if)^2 - \sum_{i,j=1}^N x_ix_j (\pp_i f)\pp_j f\\
&\ge \sum_{i=1}^N x_i (\pp_if)^2 - \sum_{i,j=1}^N x_ix_j\cdot \ff{ (\pp_i f)^2+(\pp_j f)^2} 2=\sum_{i=1}^N x_i (1-|x|_1) (\pp_i f)^2.\end{align*}
So,  $\EE_\aa^{(N)}(f,f)\le \EE_{FV}^{\aa,N}(f,f)$, and   we only need to prove the desired Nash inequality for $(\EE_\aa^{(N)},\mu_\aa^{(N)}).$  To this end, it suffices to prove
\beq\label{SP'}\mu_\aa^{(N)} (f^2)\le r\EE_\aa^{(N)}(f,f) + cr^{-p_\aa} \mu_\aa^{(N)}(|f|)^2,\ \ r\in (0,r_1],\ f\in C^1(\DD^{(N)})\end{equation}
  for some constants $c, r_1>0.$ Indeed, this inequality is equivalent to
\beq\label{SP''}\mu_\aa^{(N)} (f^2)\le r\EE_\aa^{(N)}(f,f) + c(r\land r_1)^{-p_\aa} \mu_\aa^{(N)}(|f|)^2,\ \ r>0,\ f\in C^1(\DD^{(N)}).\end{equation}
 Since by \cite[Theorem 1.1]{FMW} the generator $L_\aa^{(N)}$ has  spectral gap $\aa_{N+1}>0$, there holds
 \beq\label{SP'''} \mu_\aa^{(N)} (f^2)\le \ff 1 {\aa_{N+1}}\EE_\aa^{(N)}(f,f),\ \ r>0,\ f\in C^1(\DD^{(N)}), \mu_\aa^{(N)}(f)=0.\end{equation}
Noting that for some constant $c(r_1,\aa_{N+1})>0$ we have
$$(r\land r_1)^{-p_\aa}\le c(r_1,\aa_{N+1}) r^{-p_\aa},\ \ r\in (0, \aa_{N+1}^{-1}),$$
so that  \eqref{SP''} and \eqref{SP'''}  yield
$$\mu_\aa^{(N)} (f^2)\le r\EE_\aa^{(N)}(f,f) + c'r^{-p_\aa} \mu_\aa^{(N)}(|f|)^2,\ \ r>0,\ f\in C^1(\DD^{(N)}), \mu_\aa^{(N)}(f)=0 $$ for some constant $c'>0$. Minimizing the upper bound in $r>0$, we prove   \eqref{SB} for some constant $C>0$ and $p=p_\aa$.

To prove \eqref{SP'}  using Theorem \ref{T2.1}, we denote $x_{N+1}=1-|x|_1=1-\sum_{i=1}^N x_i$ and take
\beq\label{PHI}\phi(x)= x_{N+1}^{-1},\ \ x=(x_i)_{1\le i\le N}\in \DD^{(N)}.\end{equation} Then
\beq\label{DS}D_s:=\{\phi\le s\}= \big\{x\in \DD^{(N)}:\ x_{N+1}\ge s^{-1}\big\},\ \ s>1.\end{equation}
For the present model we have
$$\GG(\phi,\phi)(x)=\sum_{i=1}^N x_i x_{N+1} (\pp_i\phi)^2(x)=\ff{1-x_{N+1}}{x_{N+1}^3} \le s^3,\ \ x\in D_s, s>0,$$
so that 
\beq\label{DHS} h(s) :=\sup_{D_s} \GG(\phi,\phi) \le s^3,\ \ s>0.\end{equation}
To apply Theorem \ref{T2.1}, we take $$\D(\GG)=\Big\{f\in C(\DD^{(N)}; [-\infty,\infty]): f \ \text{is\ finite\ and}\  C^1 \ \text{in}\ \DD^{(N)}\setminus \{x_{N+1}=0\}\Big\},$$ and let
$$\GG(f,g)(x)= 1_{\{x_{N+1}>0\}}\sum_{i=1}^N x_i x_{N+1} (\pp_i f)(x) (\pp_i g)(x),\ \ f,g\in \D(\GG).$$
Obviously, conditions (a)-(c) hold, and the function $\phi$ in \eqref{PHI} meets the requirement of  Theorem \ref{T2.1}. In the following two subsections, we estimate $\ll(s)$ and $\bb_s$ respectively.

\subsection{Estimate on $\ll(s)$}

Let $\ll(s)= \inf\{\EE_\aa^{(N)}(f,f):\ f\in C^1(\DD^{(N)}), \mu_\aa^{(N)}(f^2)=1, f|_{D_s}=0\}$.  We will adopt the following Cheeger type  estimate $\ll(s)$. Let $$\pp D_s= \{x\in \DD^{(N)}: x_{N+1}=s^{-1}\},\ s\ge 1.$$

 \beg{lem}\label{L2.2}  If there exists a function $\psi\in C^2(\DD^{(N)}\setminus \{x_{N+1}=0\})$ such that
 \beq\label{GG} \GG(\psi,\psi)(x):= \sum_{i=1}^N x_ix_{N+1} (\pp_i \psi)^2(x)\le a_1,\ \ |L_\aa^{(N)} \psi|(x)\ge a_2, \ \ x\in D_s^c\end{equation}
 holds for some constants $a_1,a_2>0$, and that
 \beq\label{GG2} \lim_{r\to\infty} \sup_{x\in\pp D_r} x_{N+1}^{\aa_{N+1}}\sum_{i=1}^N x_i |\pp_i \psi(x)|=0,\end{equation}
 then
 $$\ll(s)\ge \ff{a_2^2}{4a_1}.$$\end{lem}

 \beg{proof} By \eqref{GG}, we assume that $L_\aa^{(N)} \psi|_{D_s^c}\ge a_2$, otherwise simply use $-\psi$  replacing $\psi$.
 Let $\si(x)={\rm diag}\big\{\ss{x_ix_{N+1}}\big\}_{1\le i\le N}$.
For any nonnegative $f\in C^1(\DD^{(N)})$ with $f|_{D_s}=0$, we have $f|_{\pp D_s}=0$, so that by integration by parts formula,
\beq\label{GG3} \beg{split} &a_2 \mu_\aa^{(N)} (f) \le \mu_\aa^{(N)} (f L_\aa^{(N)}\psi) =\lim_{r\to \infty} \int_{D_r\setminus D_s} (\rr f L_\aa^{(N)}\psi)(x)\d x\\
&\le -\mu_\aa^{(N)}\big(\<\si \nn f, \si\nn \psi\>\big) +\|f\|_\infty\limsup_{r\to\infty} \int_{\pp D_r} \sum_{i=1}^N x_ix_{N+1} \rr(x) |\pp_i \psi|(x) \d A,\end{split}\end{equation}
where $A$ is the area measure on $\pp D_r$ induced by the Lebesgue measure. We have
$$\pp D_r= \Big\{ \sum_{i=1}^N x_i= 1-r^{-1}\Big\},\ \ r\ge 2,$$ and
$$\int_{\pp D_r} \prod_{i=1}^N x_i^{\aa_i-1} \d A = (1-r^{-1})^{\sum_{i=1}^N\aa_i} \int_{\DD^{(N-1)}}\Big(1-\sum_{1\le i\le N-1}x_i\Big)^{\aa_N-1} \prod_{i=1}^{N-1} x_i^{\aa_i-1} \d x$$ is bounded in $r\ge 2$.
Combining this with \eqref{RR}, \eqref{GG}, \eqref{GG2} and \eqref{GG3}, we obtain
$$a_2 \mu_\aa^{(N)} (f)\le |\mu_\aa^{(N)}\big(\<\si \nn f, \si\nn \psi\>\big)|\le \ss{a_1 } \mu_\aa^{(N)}(|\si\nn f|).$$ Therefore, for any $f\in C^1(\DD^{(N)})$ with $f|_{D_s}=0$,
$$\mu_\aa^{(N)}(f^2)\le \ff{\ss{a_1}}{a_2} \mu_\aa^{(N)}(|\si\nn f^2|)\le \ff{2\ss{a_1}}{a_2}\ss{\mu_\aa^{(N)}(f^2)\mu_\aa^{(N)}(|\si\nn f|^2)}.$$
Noting that $\mu_\aa^{(N)}(|\si\nn f|^2)=\EE_\aa^{(N)}(f,f)$, we arrive at
$$\mu_\aa^{(N)}(f^2)\le \ff{4a_1}{a_2^2}\EE_\aa^{(N)}(f,f),\ \ f\in C_b^1(\DD^{(N)}), f|_{D_s}=0, $$
which finishes the proof. \end{proof}

\beg{lem}\label{L2.3}  There exist  constants $s_0,c_0>0$ such that
$$\ll(s)\ge c_0 s,\ \ s\ge s_0.$$\end{lem}
\beg{proof}  Let $\gg\in [\ff 1 2,1)\cap (1-\aa_{N+1},1)$. Take
$$\psi(x)= x_{N+1}^\gg, \ \ x\in \DD^{(N)}.$$ Then
\beq\label{DG} \beg{split}   \GG(\psi,\psi)(x)&=\sum_{i=1}^N x_ix_{N+1} (\pp_i \psi)^2(x)= \gg^2 (1-x_{N+1})x_{N+1}^{2\gg -1} \\
 &\le \gg^2s^{1-2\gg},\ \ x\in D_s^c,\end{split}\end{equation} and
\beq\label{DG'} \lim_{r\to\infty} \sup_{x\in D_r} x_{N+1}^{\aa_{N+1}} \sum_{i=1}^N x_i |\pp_i \psi(x)|  \le \lim_{r\to\infty} \gg s^{1-\aa_{N+1}-\gg} =0. \end{equation}
Let $s_0\ge 1$ such that
$$(1+\aa_{N+1}-\gg)  (1-s^{-1}) \ge 1-\gg + s^{-1} \sum_{i=1}^N \aa_i,\ \ s\ge s_0.$$ Then for $x\in D_s^c$ and $s\ge s_0$,
\beg{align*} L_\aa^{(N)}\psi(x)& =\sum_{i=1}^N \big\{x_i x_{N+1} \pp_i^2\psi(x) +  (\aa_ix_{N+1}-\aa_{N+1} x_i)\pp_i \psi(x)\\
&= \gg (1+\aa_{N+1}-\gg) (1-x_{N+1}) x_{N+1}^{\gg-1} -\gg x_{N+1}^\gg \sum_{i=1}^N \aa_i \\
&\ge \gg(1-\gg)s^{1-\gg}.
\end{align*}
Combining this with \eqref{DG} and \eqref{DG'}, we derive from   Lemma \ref{L2.2} that
$$\ll(s)\ge \ff{\gg^2(1-\gg)^2 s^{2(1-\gg)}}{4\gg^2s^{1-2\gg}}=  \ff{(1-\gg)^2}{4}  s,\ \ s\ge s_0.$$
 \end{proof}
\subsection{Estimate on $\bb_s(r)$}
We first present a sharp super Poincar\'e inequality for a product probability measure, then estimate $\bb_s(r)$ using a perturbation argument. Consider the following probability measures on $[0,1]$: $$\mu_i(\d s) =\aa_i s^{\aa_i-1}\d s,\ \ 1\le i\le N,$$ and let
$\mu=\prod_{i=1}^N \mu_i$ on $[0,1]^N$. We have the following result.

\beg{lem}\label{L3.1} Let $p(\aa)=\sum_{i=1}^N (\ff 1 2\lor \aa_i).$ There exists a constant $c>0$ such that
$$\mu(f^2)\le r \int_{[0,1]^N} \sum_{i=1}^N x_i (\pp_i f)^2(x) \mu(\d x) + c\big(r^{-p(\aa)} +1\big)\mu(|f|)^2,\ \ r>0, f\in C^1([0,1]^N).$$\end{lem}
\beg{proof} By Proposition \ref{PAD}, it suffices to prove that for every $1\le i\le N$ there exists a constant $c_i>0$ such that
\beq\label{SPII} \mu_i(f^2) \le r \int_0^1 s  f'(s)^2 \d s + c_i(1+r^{-(\ff 1 2\lor \aa_i)})\mu_i(|f|)^2,\ \ r>0, f\in C^1([0,1]).\end{equation}
For fixed $1\le i\le N$, we will prove this inequality using isoperimetric constants
$$\kk(r):= \inf_{I\subset [0,1], 0<\mu(I)\le r} \ff{A_i((\pp I)\setminus \{0,1\})}{\mu_i(I)},\ \ r\in (0,1/2),$$
where $A_i$ is the boundary measure induced by $\mu_i$ and the intrinsic metric of the square field $\GG_0(f,f)(s):= s (f')^2(s)$ on $[0,1]$.
Let $r\in (0, \frac{1}{2})$, for any measurable set $I \subset [0,1]$ with $\mu_i(I)=r$, we may find out $a\in (\pp I)\setminus \{0,1\}$ such that $a\ge r^{\ff 1 {\aa_i}}.$ Otherwise,    $[r^{\ff 1 {\aa_i}},1)$ is either in the interior of $I$ or in that of $I^c$. For the first case we have
$$r=\mu_i(I)\ge \aa_i \int_{r^{\ff 1 {\aa_i}}}^1 s^{\aa_i-1}\d s =1-r> \frac{1}{2}>r$$ which is a contraction; while in the second case we may find out small $\vv>0$ such that $[r^{\ff 1 {\aa_i}}-\vv, 1)\subset I^c$, so that
$$r=\mu_i(I)\le \aa_i\int_0^{r^{\ff 1 {\aa_i}}-\vv} s^{\aa_i-1}\d s = (r^{\ff 1 {\aa_i}}-\vv)^{\aa_i} <r
$$ which is again impossible.  Since the intrinsic metric induced by $\GG_0$ is
$$\d(s,t):= 2|\ss s-\ss t|,\ \ \  s,t\in [0,1],$$ the corresponding boundary measure of $\{a\}$ is given by
$$A_i(\{a\}):= \lim_{\vv\downarrow 0} \ff{\mu_i([a-\vv,a])}{2(\ss a-\ss{a-\vv})} = \ss a a^{\aa_i-1} = a^{\aa_i-\ff 1 2}.$$
Therefore,
$$\ff{A_i((\pp I)\setminus \{0,1\})}{\mu_i(I)}\ge \ff{A_i(\{a\})} r \ge r^{-(1\land \ff 1 {2\aa_i})}.$$ Hence,
$$\kk(r)\ge r^{-(1\land \ff 1 {2\aa_i})},\ \ r\in (0,1/2).$$
According to \cite[Theorem 3.4.16(1)]{Wbook}, this implies \eqref{SPII} for some constant $c_i>0$.
\end{proof}

 \beg{lem}\label{L3.2} Let $p(\aa)=\sum_{i=1}^N (\ff 1 2\lor \aa_i).$ There exist constants $c_0,s_0>0$ such that for any $s\ge s_0$,
 $$\mu_\aa^{(N)} (f^2)\le r \EE_\aa^{(N)}(f,f) + \bb_s(r) \mu_\aa^{(N)}(|f|)^2,\ \ r>0, f\in C^1(\DD^{(N)}), f|_{D_s^c}=0$$ holds for
\beq\label{LL3.2} \bb_s(r)=  c_0s^{p(\aa)+(\aa_{N+1}-1)^+} \big(r^{-p(\aa)}+ s^{p(\aa)}\big),\ \ r>0.\end{equation}
\end{lem}
\beg{proof}Let $f\in C^1(\DD^{(N)}),~f|_{D_s^c}=0$. For simplicity, we will regard $x_i$ as the function mapping $x\in \DD^{(N)}$ into $x_i, 1\le i\le N+1$. Recall that $x_{N+1}:= 1-\sum_{i=1}^N x_i$. 
Applying Lemma \ref{L3.1} to $g:= x_{N+1}^{(\aa_{N+1}-1)/2} f$ replacing $f$, which is supported on $D_s$, we may find out  constants $c_1,c_2,c_3,c_4>0$ such that for any $t>0$ and $s\ge 1$,
\beg{align*}& \mu_\aa^{(N)}(f^2)  =c_1\mu(g^2) \le c_1t \mu\Big(\sum_{i=1}^N x_i (\pp_i g)^2\Big) + c_2 (1+ t^{-p(\aa)})\mu(|g|)^2 \\
& \le t c_3 \mu_\aa^{(N)}\Big(\sum_{i=1}^N x_i\{ (\pp_i f)^2+  x_{N+1}^{-2} f^2\}\Big)+ c_2 (1+ t^{-p(\aa)})  \mu_\aa^{(N)} \big(x_{N+1}^{-\ff{\aa_{N+1}-1}2} |f|\big)^2\\
&\le c_3t s \mu_\aa^{(N)}\Big(\sum_{i=1}^N x_i x_{N+1} (\pp_i f)^2\Big)
 + c_3 t s^2 \mu_\aa^{(N)}(f^2)
 + c_2 (1+t^{-p(\aa)}) s^{(\aa_{N+1}-1)^+} \mu_\aa^{(N)}(|f|)^2\\
&\le c_3 t  s\EE_\aa^{(N)}(f,f)+ c_3t s^2 \mu_\aa^{(N)}(f^2)
 + c_2 (1+t^{-p(\aa)}) s^{(\aa_{N+1}-1)^+} \mu_\aa^{(N)}(|f|)^2.\end{align*} For any $r>0$, take
$$t= \ff{r}{2 c_3 s}\land \ff 1 {2c_3 s^2}.$$
We may find out a constant $c>0$ such that the above gives 
$$\mu_\aa^{(N)}(f^2)\le r \EE_\aa^{(N)}(f,f) + c\big(   r^{-p(\aa)}+ s^{p(\aa)}\big)s^{p(\aa)+ (\aa_{N+1}-1)^+}\mu_\aa^{(N)}(|f|)^2,\ \ r>0.$$ Therefore, the proof is finished.

\end{proof}

\subsection{Proof of \eqref{SP'}}
By \eqref{SR} and Lemma \ref{L2.3}, there exist constants $r_0,c_1>0$ such that
\beq\label{LL3.3} s_r\le  c_1 r^{-1},\ \ r\in (0,r_0].\end{equation}
Combining this with  \eqref{DHS}, we obtain
$$\ff{r h(2s_r)}{s_r^2}\le 8 r s_r \le 8c_1.$$
So, there exist constants $r_1\in (0,r_0]$ and $c_2>1$  such that for any $r\in (0,r_1]$,
$$\ff{4s_r^2+h(2s_r)r}{2s_r^2}\le c_2,\ \ \ff{r }{8+ 2rh(2s_r)s_r^{-2}} \ge \ff r {c_2}.$$
Combining these with \eqref{LL3.2} and \eqref{LL3.3},  we may find out constants $c_3,c_4>0$ such that $\bb(r)$ in \eqref{BT} satisfies
\beg{align*} \bb(r)&\le c+ 2\bb_{3s_r}(r/c_2) \le c_3 s_r^{p(\aa)+(\aa_{N+1}-1)^+}\big(r^{-p(\aa)} +s_r^{p(\aa)}\big)\\
&\le c+ c_4 r^{-\{2p(\aa)+(\aa_{N+1}-1)^+\}},\ \ r\in (0,r_1].\end{align*}
  This completes the proof since
  $$2p(\aa)+(\aa_{N+1}-1)^+= \sum_{i=1}^N 1\lor(2\aa_i) + (\aa_{N+1}-1)^+=p_\aa.$$

\section{Proof of Theorem \ref{T1.1}(2)-(3)}
\beg{proof}[Proof of Theorem \ref{T1.1}(2)] Let \eqref{SB} hold. We aim to prove $p\ge \tt p_\aa^{(1)}$ and $p\ge \tt p_\aa^{(2)}$ respectively, where
\beg{align*} &\tt p_\aa^{(1)}:= \aa_{N+1}  + \max_{1\le i\le N} \sum_{1\le j\le N, j\ne i} (1\lor \aa_j),\\
&\tt p_\aa^{(2)}:=\max_{1\le i\le N+1} \sum_{1\le j\le N+1, j\ne i} \aa_j.\end{align*}

(a) Let $1\le i_0\le N$ be such that $\aa_{i_0}= \min_{1\le i\le N}\aa_i$. Let
$$I_1=\{i_0\}\cup\{1\le i\le N: \aa_i\le 1\},\ \ I_2=\{1,\cdots, N\}\setminus I_1.$$
We have $n_1:=\# I_1\ge 1, \ \# I_2=N-n_1,$ and
\beq\label{D*0} \sum_{i\in I_2} (\aa_i-1)= \max_{1\le i\le N} \sum_{1\le j\le N, j\ne i} (\aa_j-1)^+= \max_{1\le i\le N} \sum_{1\le j\le N, j\ne i} (1\lor \aa_j) +1-N.\end{equation}
Take $h\in C^\infty(\R)$ such that $0\le h\le 1, |h'|\le 2$ and
$$h|_{(-\infty,1]}=h|_{[4,\infty)}=0,\ \ h|_{[2,3]}=1.$$
Let $\vv_N= \ff 1 {32 N^2}$ and take
\beq\label{*0D} f_\vv (x)= \bigg(\prod_{i\in I_1} h\Big(\ff{n_1^{-1}-x_i}{4N\vv}\Big)\bigg)\cdot \prod_{i\in I_2} \Big(1-\ff{x_i}{2\vv}\Big)^+,\ \ x\in \DD^{(N)}, \vv\in (0,\vv_N].\end{equation}  It is easy to see that $A_\vv:={\rm supp}f_\vv$ satisfies
\beg{align*}  A_\vv^{(1)}&:= [n_1^{-1} -12 N\vv, n_1^{-1}- 8N\vv]^{I_1}\times [\vv, 2\vv]^{I_2}\subset A_\vv\\
&\subset A_\vv^{(2)}:= [n_1^{-1}- 16 N\vv, n_1^{-1}- 4 N\vv]^{I_1}\times [0,2\vv]^{I_2}.\end{align*}
 So, for $x\in A_\vv$ we have
\beg{align*} 2N\vv &\le 1-\sum_{i\in I_1} (n_1^{-1}-4N\vv) -2  \vv (N-n_1)\\
 &\le 1-\sum_{i=1}^N x_i= x_{N+1}\le 1-\sum_{i\in I_1} (n_1^{-1} - 16 N\vv)\le 16 N^2\vv,\end{align*}  and there exist constants $c_2>c_1>0$ such that
 \beg{align*}&c_1 1_{A_\vv^{(1)}}(x)\vv^{\sum_{i\in I_2}(\aa_i-1) +\aa_{N+1}-1} \le (1_{A_\vv} \rr)(x) \le c_2 1_{A_\vv^{(2)}}(x)\vv^{\sum_{i\in I_2}(\aa_i-1) +\aa_{N+1}-1},\\
 & 1_{A_\vv^{(1)}}(x) \le f_\vv(x)\le 1_{A_\vv^{(2)}}(x),\\
 & \sum_{i=1}^N x_ix_{N+1} (\pp_if_\vv)^2(x) \le c_2\vv^{-1} 1_{A_\vv^{(2)}}(x),\ \ x\in \DD^{(N)}, \vv\in (0,\vv_N].\end{align*}
 Combining these together we may find out constants $c_3, c_4>0$ such that for any $\vv\in (0,\vv_N]$,
\beq\label{*1D} \beg{split} &\mu_\aa^{(N)}(f_\vv^2)\ge \mu_\aa^{(N)}(A_\vv^{(1)})=\int_{A_\vv^{(1)}} \rr(x) \d x \ge c_3 \vv^{\sum_{i\in I_2}(\aa_i-1) + N +\aa_{N+1}-1},\\
&\mu_\aa^{(N)}(f_\vv)^2 \le \mu_\aa^{(N)}(A_\vv^{(2)})^2\le c_4 \vv^{2N+2(\aa_{N+1}-1) +  2\sum_{i\in I_2} (\aa_i-1)},\\
&\mu_\aa^{(N)}\Big( \sum_{i=1}^N x_ix_{N+1}  (\pp_i f_\vv)^2\Big)\le c_2\vv^{-1} \mu_\aa^{(N)}(A_\vv^{(2)})  \le c_4  \vv^{\sum_{i\in I_2}(\aa_i-1) + N +\aa_{N+1}-2}.\end{split}\end{equation} Therefore, if  \eqref{SP} holds then
\beg{align*}&c_3   \vv^{\sum_{i\in I_2}(\aa_i-1) + N +\aa_{N+1}-1}\\
 &\le r  c_4 \vv^{\sum_{i\in I_2}(\aa_i-1) + N +\aa_{N+1}-2} +c_4\bb(r) \vv^{2N+2(\aa_{N+1}-1) +  2\sum_{i\in I_2} (\aa_i-1)},\ \ r>0,\vv\in (0,\vv_N].\end{align*} This is equivalent to
$$1-\ff{c_4}{c_3} r \vv^{-1} \le \ff{c_4}{c_3} \bb(r) \vv^{\sum_{i\in I_2}(\aa_i-1) + N +\aa_{N+1}-1},\ \ r>0,\vv\in (0,\vv_N].$$
Let $r_N= \ff{c_3}{2c_4} \vv_N$. For any $r\in (0, r_N]$, we take $\vv= \ff{2c_4}{c_3} r\in (0, \vv_N]$ in  the above inequality to derive from \eqref{D*0} that
$$\bb(r)\ge \ff{c_3}{2c_4}  \vv^{1-\sum_{i\in I_2}(\aa_i-1) - N -\aa_{N+1}}= c r^{-\tt p_\aa^{(1)}},\ \ r\in (0, r_N] $$ for some constant $c>0$.
Since \eqref{SB} implies \eqref{SP} for $\bb(r)= c(1+r^{-p})$ for some constant $c>0$, this implies $p\ge \tt p_\aa^{(1)}.$

(b) On the other hand,  $1\le i_0\le N+1$ be such that $\aa_{i_0}= \min_{1\le i\le N+1}\aa_i$. Let $$I=\{i: i\ne i_0, 1\le i\le N+1\}.$$  For any $\vv\in (0,1)$, we take
$$f_\vv(x)=  \prod_{i\in I} \big(\vv-x_i \big)^+,\ \ x\in \DD^{(N)}.$$ Then on the support of $f_\vv$ we have
$$x_ix_{N+1}\le \min\{x_i,x_{N+1}\}\le \vv,\ \ 1\le i\le N.$$ So, as shown in (a)  we may find out  a  constant $a\in (0,1)$ such that for all $\vv\in (0,1)$,
\beg{align*} &\mu_\aa^{(N)}(f_\vv^2)\ge a \vv^{3N+  \sum_{i\in I} (\aa_i-1) }=a \vv^{2N+ \sum_{i\in I} \aa_i},\\
&\mu_\aa^{(N)}(f_\vv)^2 \le a^{-1} \vv^{4N+   2\sum_{i\in I} (\aa_i-1)}=a^{-1} \vv^{2N+ 2 \sum_{i\in I} \aa_i},\\
&\mu_\aa^{(N)}( |\si \nn f_\vv|^2)\le a^{-1}  \vv^{3N+  \sum_{i\in I} (\aa_i-1) -1} =a^{-1} \vv^{2N+  \sum_{i\in I} \aa_i  -1},\end{align*} and these together with  \eqref{SB} imply
$$ p \ge  \sum_{i\in I} \aa_i= \max_{1\le i\le N+1} \sum_{j\ne i, 1\le j\le N+1}\aa_j=\tt p_\aa^{(2)}.$$\end{proof}

\beg{proof}[Proof of Theorem \ref{T1.1}(3)] Let $f_\vv$ be in \eqref{*0D}. We have 
$$\EE_{FV}^{\aa,N}(f_\vv,f_\vv)\le \mu_\aa^{(N)}\bigg(\sum_{i=1}^N x_i (\pp_i f_\vv)^2\bigg) \le c_4 \vv^{\sum_{i\in I_2} (\aa_i-1)+N+\aa_{N+1}-1}$$ for some constant $c_4>0$. Combining this with the first two lines in \eqref{*1D}, we derive from \eqref{SP} with $\EE_{FV}^{\aa,N}$ replacing $\EE_\aa^{(N)}$ that 
$$c_3\vv^{\sum_{i\in I_2} (\aa_i-1) +N +\aa_{N+1}-1} \le c_4 \vv^{\sum_{i\in I_2} (\aa_i-1) +N +\aa_{N+1}-3}
+c_4 \vv^{2\sum_{i\in I_2} (\aa_i-1) +2N +2\aa_{N+1}-2}\bb(r),$$ thus,
$$ 1-\ff{c_4r}{c_3 \vv^2}\le \ff{c_4}{c_3}\bb(r) \vv^{\sum_{i\in I_2} (\aa_i-1) +N +\aa_{N_1}-1}.$$
Taking $\vv= \big(\ff{2c_3r}{c_4}\big)^{\ff 1 2}$ for small $r>0$ we arrive at
$$\bb(r)\ge cr^{-\ff 1 2 (N+\aa_{N+1} -1 +\sum_{i\in I_2} (\aa_i-1))}.$$
Combining this with \eqref{SB} implies $p\ge \ff 1 2\aa_{N+1} +\ff 1 2 \sum_{i=1}^N (1\lor \aa_i).$

On the other hand, take
$$f_\vv(x)= \prod_{1\le i\le N} (\vv-x_i)^+$$ for small $\vv>0$. Then there exists a constant $a\in (0,1)$ such that for small $\vv>0$ we have \beg{align*} &\mu_\aa^{(N)}(f_\vv^2)\ge  a \vv^{2N+ \sum_{1\le i\le N} \aa_i},\\
&\mu_\aa^{(N)}(f_\vv)^2 \le a^{-1} \vv^{2N+ 2 \sum_{1\le i\le N} \aa_i},\\
&\EE_{FV}^{\aa,N}(f_\vv,f_\vv)\le \mu_\aa^{(N)}\bigg( \sum_{i=1}^N x_i (\pp_i f_\vv)^2\bigg)\le a^{-1}  \vv^{2N+  \sum_{1\le i\le N} \aa_i-1},\end{align*} so that \eqref{SB} for $\EE_{FV}^{\aa,N}$ implies
$p \ge  \sum_{1\le i\le N} \aa_i.$\end{proof}

\beg{thebibliography}{99}

\bibitem{Bakry} D. Bakry, I. Gentil, M. Ledoux, \emph{Analysis and Geometry of Markov Diffusion Operators}, Springer 2014.

\bibitem{BR} J. Bakosi, J. R. Ristorcelli, \emph{A stochastic diffusion process for the Dirichlet distribution,} International J. Stoch. Anal. 2013, Article ID 842981, 7 pages.



\bibitem{ConMoi69} R. J. Connor,  J. E. Mosimann,
 \emph{Concepts of independence for proportions with a generalization of the Dirichlet distribution,}  J. Amer. Statist. Assoc. 64(1969), 194-206.

 \bibitem{Davies} E. B. Davies, B. Simon, \emph{Ultracontractivity and the heat kernel
for Schr\"odinger operators and Dirichlet Laplacians,} J. Funct. Anal. 59(1984), 335--395.

\bibitem{QQ} C. L. Epstein, R. Mazzeo, \emph{ Wright-Fisher diffusion in one dimension,}  SIAM J. Math. Anal.   42(2010),  568--608.




\bibitem{FMW} S. Feng, L. Miclo, F.-Y. Wang, \emph{ Poincare inequality for Dirichlet distributions and infinite-dimensional generalizations,} Lat. Am. J. Probab. Math. Stat. 14(2017), 361--380.

 \bibitem{FW07}  S. Feng, F.-Y. Wang,  \emph{A class of infinite-dimensional diffusion processes with connection to population genetics,} J. Appl. Probab. 44(2007), 938--949.

 \bibitem{FW14}  S. Feng,   F.-Y. Wang,   \emph{Harnack inequality and applications for infinite-dimensional GEM   processes,} Potential Anal.
  44(2016),  137--153


 \bibitem{Gross} L. Gross, \emph{ Logarithmic Sobolev inequalities
and contractivity properties of semigroups,}  in
$``$Dirichlet Forms'',  Lecture Notes in
Math.  1563 (Springer, Berlin), pp. 54--88.

\bibitem{Jac01} M. Jacobsen, \emph{Examples of multivariate diffusions: time-reversibility; a Cox-Ingersoll-Ross type process,} Department of Theoretical Statistics, Preprint 6, University of Copenhagen, 2001.


   \bibitem{P1} N. L. Johnson, \emph{An approximation to the multinomial distribution, some properties and applications,} Biometrika, 47(1960), 93--102.

   \bibitem{Mi1} L. Miclo, \emph{About projections of logarithmic Sobolev inequalities,}  Lecture Notes in Math. 1801 (J. Az\'ema, M. \'Emery, M. Ledoux, M. Yor Eds), pp. 201--221, 2003, Springer.

 \bibitem{Mi2} L. Miclo, \emph{Sur l'in\'egalit\'e de Sobolev logarithmique des op\'erateurs de Laguerre \`{a} petit param\'etre, } Lecture Notes in Math. 1801 (J. Az\'ema, M. \'Emery, M. Ledoux, M. Yor Eds), pp. 222--229, 2003, Springer.

 \bibitem{P2} J. E. Mosimann, \emph{On the compound multinomial distribution, the multivariate-distribution, and correlations among proportions,}
 Biometrika, 49(1962), 65--82.

 \bibitem{SN} N. Shimakura, \emph{Equations diff\'erentielles provenant de la g\'enetique des populations,} T\^ohoka Math. J. 29(1977), 287--318.

  \bibitem{S} W. Stannat, \emph{On validity of the log-Sobolev
      inequality for symmetric Fleming-Viot operators,}  Ann.
      Probab.
      28(2000), 667--684.
\bibitem{W00a}
F.-Y. Wang, \emph{ Functional inequalities for empty essential
spectrum, }  J. Funct. Anal. 170(2000), 219--245.

\bibitem{W00b} F.-Y. Wang, \emph{ Functional inequalities, semigroup properties
and spectrum estimates,}   Infin. Dimens. Anal. Quant. Probab.
Relat. Topics 3(2000), 263--295.

\bibitem{Wbook} F.-Y. Wang,   \emph{Functional Inequalities, Markov Semigroups and Spectral Theory,}   Science Press 2005.



\end{thebibliography}
\end{document}